\documentclass[twoside]{aiml22}

\usepackage{aiml22macro}

\usepackage{graphicx}
\usepackage{amsmath}
\usepackage{amssymb}

\usepackage{stackengine} 
\usepackage{multicol}

\newcommand\loz{\stackMath\mathbin{\stackinset{c}{0ex}{c}{0ex}{\ast}{\lozenge}}}

\newcommand{\ignore}[1]{}

\newcommand{\sfont}[1]{\mathfrak{#1}}

\newcommand{\dn}{\square}

\newcommand{\val}[1]{\llbracket #1 \rrbracket}
\usepackage{tikz}
\usepackage{float}
\usepackage{tikz-cd}
\usepackage{caption}
\usepackage[font={small,it}]{caption}
\usepackage{siunitx}
\usepackage{centernot}
\usepackage{stmaryrd}
\usepackage{hyperref}

\usepackage{wasysym}

\NewDocumentCommand{\bc}{}{{\text{\normalfont\newmoon}}}
\newcommand{\fulllan}{\mathcal{L}^{\bc\blacksquare}_{\square}}

\newcommand{\lb}{\langle}
\newcommand{\rb}{\rangle}

\def\lastname{Montacute}

\begin{document}

\begin{frontmatter}
  \title{Cantor Derivative Logic in Topological Dynamics}
 % \author{David Fern\'andez-Duque}
 % \address{Affiliation \\ Address \\ Address}
  \author{Yo\`av Montacute}
  \address{University of Cambridge \\ Cambridge \\ United Kingdom}

  \begin{abstract}
   Topological semantics for modal logic based on the Cantor derivative operator gives rise to derivative logics, also referred to as $d$-logics.
Unlike logics based on the topological closure operator, $d$-logics have not previously been studied in the framework of dynamic topological systems (DTSs), which are pairs $\langle X,f \rangle$ consisting of a topological space $X$ equipped with a continuous function $f\colon X\to X$.
We introduce the logics $\bf{wK4C}$, $\bf{K4C}$ and $\bf{GLC}$ and show that they all have the finite Kripke model property and are sound and complete with respect to the $d$-semantics in this dynamical setting. We also prove a general result for the case where $f$ is a homeomorphism, which yields soundness and completeness for the corresponding systems $\bf{wK4H}$, $\bf{K4H}$ and $\bf{GLH}$.
Of special interest is $\bf{GLC}$, which is the $d$-logic of all DTSs based on a scattered space.  We use the completeness of $\bf{GLC}$ and the properties of scattered spaces to demonstrate the first sound and complete dynamic topological logic in the original trimodal language. In particular, we show that the version of $\mathbf{DTL}$ based on the class of scattered spaces is finitely axiomatisable over the original language, and that the natural axiomatisation is  sound and complete. 

  \end{abstract}

 \end{frontmatter}

\section{Introduction}

\emph{Dynamic topological systems} (DTSs) are mathematical models of processes that may be iterated indefinitely.
Formally, they are defined as pairs $\lb\mathfrak X,f\rb$ consisting of a topological space $\mathfrak{X}=\langle X,\tau\rangle$ and a continuous function $f\colon X\to X$; the intuition is that points in the space $\mathfrak X$ `move' along their orbit, $x,f(x),f^2(x),\ldots$ which usually simulates changes in time.
{\em Dynamic topological logic} ($\bf DTL$) combines modal logic and its topological semantics with linear temporal logic (see Pnueli \cite{ltl}) in order to reason about dynamical systems.

Due to their rather broad definition, dynamical systems are routinely used in many pure and applied sciences, including computer science.
Such applications raise a need for effective formal reasoning about topological dynamics.
Here, we may take a cue from modal logic and its topological semantics.
The study of the latter dates back to McKinsey and Tarski~\cite{Tarski}, who proved that the modal logic $\bf S4$ is complete for a wide class of spaces, including the real line.
Artemov, Davoren and Nerode~\cite{artemov} extended $\bf S4$ with a `next' operator in the spirit of $\bf LTL$, producing the logic $\mathbf{S4C}$. 
They proved that this logic is sound and complete with respect to the class of all DTSs.
The system $\mathbf{S4C}$ was enriched with the `henceforth' operator by Kremer and Mints \cite{kremer}, who named the new logic {\em dynamic topological logic} ($\mathbf{DTL}$). 
Later, Konev et al.\ \cite{konev} showed that $\mathbf{DTL}$ is undecidable, and Fernández-Duque~\cite{david} showed that it is not finitely axiomatisable on the class of all DTSs.

The aforementioned work on dynamic topological logic interprets the modal operator $\lozenge$ as a closure operator.
However, McKinsey and Tarski had already contemplated semantics that are instead based on the Cantor derivative \cite{Tarski}. The {\em Cantor derivative} of a set $A$, usually denoted by $d(A)$, is the set of points $x$ such that $x$ is in the closure of $A\setminus \{x\}$.
This interpretation is often called {\em $d$-semantics} and the resulting logics are called {\em $d$-logics.}
These logics were first studied in detail by Esakia, who showed that the $d$-logic $\mathbf{wK4}$ is sound and complete with respect to the class of all topological spaces \cite{Esakia2}.
It is well-known that semantics based on the Cantor derivative are more expressive than semantics based on the topological closure.
For example, consider the property of a space $\mathfrak X$ being {\em dense-in-itself}, meaning that $\mathfrak X$ has no isolated points.
The property of being dense-in-itself cannot be expressed in terms of the closure operator, but it {\em can} be expressed in topological $d$-semantics by the formula $\lozenge\top$.

Logics based on the Cantor derivative appear to be a natural choice for reasoning about dynamical systems. However, there are no established results of completeness for such logics in the setting of these systems.
Our goal is to prove the finite Kripke model property, completeness and decidability of logics with the Cantor derivative operator and the `next' operator $\bc$ over some prominent classes of DTSs. Namely, those based on arbitrary spaces, $T_D$ spaces (spaces validating the $4$ axiom $\square p \to \square\square p$) and {\em scattered spaces}.
Scattered spaces are topological spaces where every non-empty subspace has an isolated point.
The reason for considering scattered spaces is to circumvent the lack of finite axiomatisability of $\mathbf{DTL}$ by restricting to a suitable subclass of all DTSs.
In the study of dynamical systems and topological modal logic, one often works with {\em dense-in-themselves} spaces.
This is a sensible consideration when modelling physical spaces, as Euclidean spaces are dense-in-themselves.
However, some technical issues that arise when studying $\mathbf{DTL}$ over the class of all spaces disappear when restricting our attention to scattered spaces, which in contrast have many isolated points.
Further, we consider DTSs where $f$ is a {\em homeomorphism,} i.e.\ where $f^{-1}$ is also a continuous function.
Such DTSs are called {\em invertible}. 

The basic dynamic $d$-logic we consider is $\mathbf{wK4C}$, which consists of $\mathbf{wK4}$ and the temporal axioms for the continuous function $f$.
In addition, we investigate two extensions of $\mathbf{wK4C}$:  $\mathbf{K4C}$ and $\mathbf{GLC}$.
As we will see, $\mathbf{K4C}$ is the $d$-logic of all DTSs based on a $T_D$ space, and $\mathbf{GLC}$ is the $d$-logic of all DTSs based on a scattered space.
Scattered spaces have gathered attention lately in the context of computational logic, as they may be used to model provability in formal theories~\cite{abashidze1985}, leading to applications in characterising their provably total computable functions \cite{Beklemishev04}. Modal logic on scattered spaces enjoys definable fixed points~\cite{SambinV82}, relating it to the topological $\mu$-calculus~\cite{BBFMu}.
The latter is particularly relevant to us, as the expressive power gained by topological fixed points, including the tangled operators of $\mathbf{DTL}^*$, is absent in this setting.
As the logic of scattered spaces is the G\"odel-L\"ob modal logic $\bf GL$, we refer to the dynamic topological logic of scattered spaces as {\em dynamic G\"odel-L\"ob logic} ($\bf DGL$).

Our goal is to demonstrate that the standard finite axiomatisation of $\bf DGL$ is sound and complete, leading to the first complete trimodal dynamic topological logic, as well as the first such logic combining the Cantor derivative with the infinitary `henceforth' operator from $\bf LTL$.
By the `standard axioms' we refer to the combination of the well-known axiomatisation of $\bf GL$ with $\bf LTL$ axioms for the temporal operators and $( \newmoon p \wedge \newmoon \Box p) \to\Box\newmoon p $ -- a variant of the continuity axiom of Artemov et al.~adapted for the Cantor derivative.
The proof of completeness employs various advanced techniques from modal logic, including an application of Kruskal's theorem in the spirit of the work of Gabelaia et al.\ \cite{pml}.

This paper summarises the results of Yoàv Montacute's Master's Thesis \cite{montacutethesis} and the consequent results it yielded \cite{D-FM21,untangled}.

\section{Dynamic topological logic with the Cantor derivative}

Given a non-empty set $\mathsf{PV}$ of propositional variables, the language $\mathcal{L}_{\square}^{\bc}$ is defined recursively as follows:
$$
\varphi::= p \; | \; \varphi\wedge \varphi \; | \; \neg\varphi \; | \; \square\varphi  \; | \; \bc\varphi \; | \; \blacksquare\varphi,$$
%$$p \in \mathsf{PV}.$$
where $p \in \mathsf{PV}.$
It consists of the Boolean connectives $\wedge$ and $\neg$, the temporal modalities `next' $\newmoon$ and `henceforth' $\blacksquare$ with its dual `eventually' $\blacklozenge:=\neg\blacksquare\neg$, and the spatial modality $\square$ for the co-derivative with its dual the Cantor derivative $\lozenge:=\neg\square\neg$. We define other connectives (e.g.\; $\vee$, $\rightarrow$) in the usual way.

\begin{definition}
	A \emph{dynamic topological model} is a tuple $\sfont M=\langle X,\tau,f,\nu\rangle$, where $\langle X,\tau,f\rangle$ is a DTS and $\nu:\mathsf{PV}\rightarrow\wp(X)$ is a {\em valuation function}. 
	Given $\varphi\in\fulllan$, we define the \emph{truth set} $\val\varphi \subseteq X$ of a formula $\varphi$ as follows:
	\setlength\columnsep{-3mm}
	\begin{multicols}{2}
\noindent
	\begin{itemize}
		\item $\val p=\nu(p)$;
		\item $\val{\neg\varphi} = X\backslash \val\varphi $;
		\item $\val{\varphi\wedge\psi}=\val\varphi \cap \val \psi $;
		\item $\val{\square\varphi}=d(\val{\neg\varphi})^\complement$;
		\item $\val{\newmoon\varphi} =f^{-1}(\val\varphi )$;
		\item $\val{\blacksquare\varphi} =\bigcap_{n\geq 0}f^{-n}(\val\varphi)$.
	\end{itemize}	
	\end{multicols}

\end{definition}
Let us list the axiom schemata and rules that we will consider in this paper:
\begin{multicols}2
\begin{description}
\item{\rm Taut} $:= \text{All propositional tautologies}$
\item{\rm K} $:= \dn(\varphi\to\psi)\to(\dn\varphi\to \dn\psi)$
\item{\rm T} $:=    \dn\varphi \to \varphi $
\item{\rm w4} $ :=  \varphi \wedge  \dn\varphi \to\dn\dn\varphi $
\item{\rm L} $:=    \square(\square \varphi \rightarrow \varphi) \rightarrow \square \varphi  $
\item{\rm 4} $  :=    \dn\varphi \to\dn\dn\varphi $
\item{${\rm Next}_\neg$} $:=\neg\bc\varphi\leftrightarrow\bc\neg\varphi$
\item{${\rm Next}_\wedge$} $:=\bc (\varphi\wedge\psi)\leftrightarrow \bc \varphi \wedge\bc \psi $
\item{\rm C} $:= \bc\varphi\wedge \bc\dn\varphi\to  \dn\bc\varphi$
\item{\rm H} $ := \dn\bc\varphi\leftrightarrow\bc\dn\varphi$
\item{\rm MP} $:= \dfrac{\varphi \ \ \varphi\to \psi}\psi$
\item{${\rm Nec}_\dn$} $:= \dfrac{\varphi }{\dn \varphi}$
\item{${\rm Nec}_\bc$} $:= \dfrac{\varphi }{\bc \varphi}$
\item{${\rm K}_\blacksquare$} $:=\blacksquare(\varphi\to\psi)\to (\blacksquare\varphi\to\blacksquare\psi)$
\item{${\rm Fix}_\blacksquare$} $:=\blacksquare\varphi\rightarrow(\varphi\wedge\newmoon\blacksquare\varphi)$
\item{${\rm Ind}_\blacksquare$} $:=\blacksquare(\varphi\rightarrow\newmoon\varphi)\rightarrow(\varphi\rightarrow\blacksquare\varphi)$
\end{description}
\end{multicols}

We define $\mathbf{wK4}:= \mathbf{K}+{\rm w4}$, $\mathbf{K4}:= \mathbf{K}+{\rm 4}$,
$\mathbf{S4}:=  \mathbf{K4}+{\rm T}$ and $\mathbf{GL}:=  \mathbf{K4}+{\rm L}$.
These are well known logics over $\mathcal L_\square$ which characterise different classes of topological spaces.
In addition, for a logic $\Lambda$ over $\mathcal L_\square$, $\Lambda\mathbf{F}$ is the logic over $\mathcal L^\bc_\square$ given by
$\Lambda\mathbf{F} := \Lambda+{\rm Next}_\neg+{\rm Next}_\wedge+{\rm Nec}_\bc$.
This simply adds axioms from linear temporal logic to $\Lambda$, which hold whenever $\bc$ is interpreted using a function.
We define $\Lambda{\bf C}:=\Lambda\mathbf{F}+\rm C$ and $\Lambda{\bf H}:=\Lambda\mathbf{F}+\rm H$, which as we will see correspond to topological spaces with a continuous function or a homeomorphism respectively.
Finally, let 
$ \bf{DGL}:=\bf{GLC}+{\rm K}_\blacksquare+{\rm Fix}_\blacksquare+{\rm Ind}_\blacksquare$.

The logic $\bold{K4}$ includes the axiom $\square p \to \square\square p$, which is not valid over the class of all topological spaces.
The class of spaces satisfying this axiom is denoted by $T_D$ and defined as the class of spaces in which every singleton is the result of an intersection between an open set and a closed set.
Moreover, Esakia showed that this is the $d$-logic of transitive Kripke frames \cite{EsakiaAlgebra}.

Many familiar topological spaces, including Euclidean spaces, satisfy the $T_D$ property, making $\mathbf{K4}$ central in the study of topological modal logic.
A somewhat more unusual class of spaces, which is nevertheless of particular interest to us, is the class of {\em scattered spaces.}

\begin{definition}
A topological space $\lb X,\tau\rb$ is {\em scattered} if for every $S\subseteq X$,
$ S\subseteq d(S) \text{ implies } S=\varnothing.$
\end{definition}
This is equivalent to the more common definition of a scattered space where a topological space is called scattered if every non-empty subset has an isolated point.
Scattered spaces are closely related to converse well-founded relations.

\begin{lemma}\label{lemScatteredKrip}
If $\langle W,\sqsubset \rangle$ is an irreflexive Kripke frame, then $\langle W,\tau_\sqsubset \rangle$ is scattered iff $\sqsubset $ is converse well-founded.
\end{lemma}

\begin{theorem}[Simmons \cite{simmons} and Esakia \cite{esakia}]\label{thmGLComp}
$\mathbf{GL}$ is the $d$-logic of all scattered topological spaces, as well as the $d$-logic of all converse well-founded Kripke frames and the $d$-logic of all finite, transitive, irreflexive Kripke frames.
\end{theorem}

Aside from its topological interpretation, the logic $\mathbf{GL}$ is especially interesting since it is also the logic of provability in Peano arithmetic, as was shown by Solovay \cite{solovay1976} (see also \cite{boolo}).
Meanwhile, logics with the $\rm C$ and $\rm H$ axioms correspond to classes of DTSs.

\clearpage 

\begin{theorem}[\cite{montacutethesis,D-FM21}]\

\begin{enumerate}

\item $\bf wK4C$ is sound and complete for
\begin{enumerate}

\item the class of finite dynamic $\bf wK4$ frames.

\item the class of finite DTSs, with respect to $d$-semantics.

\end{enumerate}

\item $\bf K4C$ is sound and complete for
\begin{enumerate}

\item the class of finite dynamic $\bf K4$ frames.

\item the class of DTSs based on a $T_D$ space, with respect to $d$-semantics.

\end{enumerate}

\item $\bf GLC$ is sound and complete for
\begin{enumerate}

\item the class of finite dynamic $\bf GL$ frames.

\item the class of finite DTSs based on a scattered space, with respect to $d$-semantics.

\end{enumerate}

\end{enumerate}
Furthermore, the same results analogously hold when the logics are replaced with $\mathbf{wK4H}$, $\mathbf{K4H}$ and $\mathbf{GLH}$, and the  Kripke frames and topological spaces are replaced with their invertible counterparts.
\end{theorem}

Kremer and Mints \cite{kremer} suggested adding the `henceforth' operator, $\blacksquare$, from Pnueli's linear temporal logic ($\bf LTL$) \cite{ltl}, leading to a trimodal system $\bf DTL$.
They offered an axiomatisation for $\bf DTL$, but Fern\'andez-Duque proved that it is incomplete; in fact, $\bf DTL$ is not finitely axiomatisable \cite{david}.
Fern\'andez-Duque also showed that $\bf DTL$ enjoys a natural axiomatisation when extended with the {\em tangled closure} \cite{david3}.

\begin{definition}
Let $\lb X,\tau\rb$ be a topological space and let $\mathcal{S}\subseteq \wp(X)$. Given $A\subseteq X$, we say that $\mathcal{S}$ is tangled in $A$ if for all $S\in \mathcal{S}$,
$A\subseteq d(A\cap S)$.
We define the \emph{tangled derivative} of $\mathcal{S}$ as 
$$ \mathcal{S}^*:=\bigcup \{ A\subseteq X : \mathcal{S}\text{ is tangled in }A\}.$$ 
\end{definition}

Fern\'andez-Duque's axiomatisation is based on the extended language with the tangled operator $\loz$.

\begin{definition}

For every model $\mathfrak M$,
$\| {\loz\{\varphi_1,\ldots,\varphi_n\}}\| = \{\|\varphi_1\|,\ldots,\|\varphi_n\|\}^* .$
\end{definition}

 Unlike the complete axiomatisation of $\mathbf{DTL}$ that requires the tangled operator, in the case of $\bold{DGL}$, we are able to avoid this and use the original spatial operators alone. This is due to the following:

\begin{theorem}
Let $\mathfrak{X}=\lb X,\tau \rb$ be a scattered space and let $\{\varphi_1,\dots,\varphi_n\}$ be a set of formulas. Then
$$ \loz\{\varphi_1,\dots,\varphi_n\}\equiv \bot.$$
\end{theorem}

Given this result together with the finite model property of $\bf{GLC}$, we can use an adapted version of the axiomatic system of Kremer and Mints \cite{kremer} in order to provide a finite axiomatisation for $\mathbf{DGL}$. 
\begin{theorem}[\cite{untangled}]
A formula $\varphi\in\mathcal \fulllan$ is valid on the class of DTSs based on a scattered space iff it is derivable in $\mathbf{DGL}$.
\end{theorem}
%% Bibliography
%% Make sure to use the bibliographystyle aiml22.
\bibliographystyle{aiml22}
\bibliography{aiml22}

\begin{thebibliography}{10}
\expandafter\ifx\csname url\endcsname\relax
  \def\url#1{\texttt{#1}}\fi
\expandafter\ifx\csname urlprefix\endcsname\relax\def\urlprefix{URL }\fi
\newcommand{\enquote}[1]{``#1''}

\bibitem{abashidze1985}
Abashidze, M., \emph{Ordinal completeness of the {G}\"odel-{L}\"ob modal
  system}, Intensional Logics and the Logical Structure of Theories  (1985),
  pp.~49--73, in Russian.

\bibitem{artemov}
Artemov, S., J.~Davoren and A.~Nerode, \emph{Modal logics and topological
  semantics for hybrid systems}, Technical Report MSI 97-05  (1997).

\bibitem{BBFMu}
Baltag, A., N.~Bezhanishvili and D.~Fern{\'{a}}ndez{-}Duque, \emph{The
  topological mu-calculus: completeness and decidability}, in: \emph{36th
  Annual {ACM/IEEE} Symposium on Logic in Computer Science, {LICS} 2021, Rome,
  Italy, June 29 - July 2, 2021} (2021), pp. 1--13.

\bibitem{Beklemishev04}
Beklemishev, L.~D., \emph{Provability algebras and proof-theoretic ordinals,
  {I}}, Ann. Pure Appl. Log. \textbf{128} (2004), pp.~103--123.

\bibitem{boolo}
Boolos, G., \emph{The logic of provability}, The American Mathematical Monthly
  \textbf{91} (1984), pp.~470--480.
\newline\urlprefix\url{https://doi.org/10.1080/00029890.1984.11971467}

\bibitem{esakia}
Esakia, L., \emph{Diagonal constructions, {L}{\"o}b’s formula and
  {C}antor’s scattered spaces}, Studies in logic and semantics \textbf{132}
  (1981), pp.~128--143.

\bibitem{Esakia2}
Esakia, L., \emph{Weak transitivity--restitution}, in: Nauka, editor,
  \emph{Study in Logic}, 2001, pp. 244--254, in Russian.

\bibitem{EsakiaAlgebra}
Esakia, L., \emph{Intuitionistic logic and modality via topology}, Ann. Pure
  Appl. Log. \textbf{127} (2004), pp.~155--170.

\bibitem{david}
Fern{\'{a}}ndez{-}Duque, D., \emph{Non-finite axiomatizability of dynamic
  topological logic}, {ACM} Transactions on Computational Logic \textbf{15}
  (2014), pp.~4:1--4:18.

\bibitem{D-FM21}
Fern\'{a}ndez-Duque, D. and Y.~Montacute, \emph{{Dynamic Cantor Derivative
  Logic}}, in: F.~Manea and A.~Simpson, editors, \emph{30th EACSL Annual
  Conference on Computer Science Logic (CSL 2022)},  Leibniz International
  Proceedings in Informatics (LIPIcs)  \textbf{216} (2022), pp. 19:1--19:17.
\newline\urlprefix\url{https://drops.dagstuhl.de/opus/volltexte/2022/15739}

\bibitem{david3}
Fernández-Duque, D., \emph{A sound and complete axiomatization for dynamic
  topological logic}, The Journal of Symbolic Logic \textbf{77} (2012),
  p.~947–969.

\bibitem{untangled}
Fernández-Duque, D. and Y.~Montacute, \emph{Untangled: A complete dynamic
  topological logic} (2022).
\newline\urlprefix\url{https://arxiv.org/abs/2204.08374}

\bibitem{pml}
Gabelaia, D., A.~Kurucz, F.~Wolter and M.~Zakharyaschev, \emph{Non-primitive
  recursive decidability of products of modal logics with expanding domains},
  Annals of Pure and Applied Logic \textbf{142} (2006), pp.~245--268.

\bibitem{konev}
Konev, B., R.~Kontchakov, F.~Wolter and M.~Zakharyaschev, \emph{Dynamic
  topological logics over spaces with continuous functions}, in:
  G.~Governatori, I.~Hodkinson and Y.~Venema, editors, \emph{Advances in Modal
  Logic} (2006), pp. 299--318.

\bibitem{kremer}
Kremer, P. and G.~Mints, \emph{Dynamic topological logic}, Annals of Pure and
  Applied Logic \textbf{131} (2005), pp.~133--158.

\bibitem{Tarski}
McKinsey, J. and A.~Tarski, \emph{The algebra of topology}, Annals of
  Mathematics \textbf{2} (1944), pp.~141--191.

\bibitem{montacutethesis}
Montacute, Y., \enquote{Chaos and Cantor derivative logic in topological
  dynamics,} Master's thesis.

\bibitem{ltl}
Pnueli, A., \emph{The temporal logic of programs}, in: \emph{Proceedings 18th
  {I}{E}{E}{E} Symposium on the Foundations of {C}{S}}, 1977, pp. 46--57.

\bibitem{SambinV82}
Sambin, G. and S.~Valentini, \emph{The modal logic of provability. the
  sequential approach}, J. Philos. Log. \textbf{11} (1982), pp.~311--342.

\bibitem{simmons}
Simmons, H., \emph{Topological aspects of suitable theories}, Proceedings of
  the Edinburgh Mathematical Society \textbf{19} (1975), p.~383–391.

\bibitem{solovay1976}
Solovay, R.~M., \emph{Provability interpretations of modal logic}, Israel
  journal of mathematics \textbf{25} (1976), pp.~287--304.

\end{thebibliography}

\end{document}